\documentclass{amsart}
\usepackage{amssymb}
\newcommand{\N}{\mathbb{N}}
\newcommand{\Z}{\mathbb{Z}}
\newcommand{\C}{\mathbb{C}}
\newcommand{\R}{\mathbb{R}}
\newcommand{\Q}{\mathbb{Q}}

\renewcommand{\P}{\mathbb {P}}

\newcommand{\G}{\mathcal{G}}

\newcommand{\res}{\mathrm{res}}
\newtheorem{Theo}{Theorem}
\newtheorem{conj}{Conjecture}

\newtheorem{Cor}{Corollary}
\newtheorem{Prop}{Proposition}

\begin{document}
\title{Analytic Continuation of some Zeta Functions}
\author[G. Bhowmik]{Gautami Bhowmik\ (Kanji name here)}
\maketitle
\section{Introduction}

The contents of this paper  
were presented as lectures at the Miura Winter School on Zeta and
$L$-functions held in 2008. 
Though the analytic continuation of zeta functions beyond its region of
absolute convergence is a fundamental question, in general not much is known
about the conditions that guarantee a meromorphic continuation. It is also
interesting to know how far such a function can be continued, that is where
the natural boundary of analytic continuation lies.

The choice of functions that are
 considered here are  `arbitrary', that is a matter
of personal taste and expertise. Most of the work reported  is on what I have studied or actually
contributed to together with my co-authors. The word  `some' in the title is to
indicate that though the paper is expository, it is not
exhaustive.  Only outlines of proofs have sometimes been provided.

In the first part we consider Euler products.
One of the most
important applications of zeta functions is
the asymptotic estimation of the sum of its coefficients via Perron's
formula, that is, the use of the equation 
\[
\sum_{n\leq x} a_n = \frac{1}{2\pi i}\int\limits_{c-i\infty}^{c+i\infty}
\Big(\sum_{n\geq 1}\frac{a_n}{n^s}\Big)\frac{x^s}{s}\;ds.
\]
To use this
relation, one usually shifts the path of integration to the left,
thereby reducing the contribution of the term $x^s$. This becomes
possible only if the function $D(s)=\sum\frac{a_n}{n^s}$ is holomorphic on the
new path. In Section 3  details of certain examples from height 
zeta functions and zeta
functions of groups have been given.

Clearly all zeta functions do not have Euler product expansions,
one important class of examples being  multiple zeta functions 
which have been studied 
often in recent years. Not many general methods exist and here I treat the
case of the Goldbach generating function  associated to $G_r(n)$, the number of
representations of $n$ as the sum of $r$ primes
\[
\sum_{k_1=1}^\infty\dots\sum_{k_r=1}^ \infty 
\frac{\Lambda(k_1)\dots\Lambda(k_r)}{(k_1+k_2+\dots+k_r)^s} 
 = \sum_{n=1}^\infty\frac{G_r(n)}{n^s}.
\]
where $\Lambda$ is the classical von-Mangoldt function.

In almost all examples the natural boundary, if it can be obtained,
corresponds to the intuitively expected boundary and this can in fact be proved in
a probabilistic sense. However one of the
difficulties in actually obtaining the boundary is that
our analyses often depend on the
distribution of zeros of the
Riemann zeta function, and thus on yet unproved hypotheses (see, for example, 
Theorem~\ref{Thm:SpecialBoundary} or Theorem~\ref{contgold} below).

I would like to thank Jean-Pierre Kahane for his comments on
Theorem~\ref{thm:RandomCont} and to 
Kohji Matsumoto for honouring me with a kanji
name. Qu'ils soient ici remerci\'es !

\section{Euler products}
 
Many Dirichlet-series occurring in practice satisfy an Euler product
and if this product is simple we often get some information on the domain
of convergence of the Dirichlet series. Among such cases is the product over all primes $p$ of a
polynomial in $p^{-s}$. One of the
oldest ideas 
is due to Estermann  \cite{Est} who
obtained a precise criterion for the continuation to the whole complex plane
of the
Euler product of an integer polynomial in $p^{-s}$. He proved the existence of the
following dichotomy : 
\begin {Theo}\label{est} Let $$ h(X)=1+a_1X+\cdots+a_dX^d =\prod_{j=1}^d(1-\alpha_jX) \in \Z[X]$$
then $Z(h;s)=\prod _p h(p^{-s})$ is absolutely convergent for $\Re(s)
>1$ and can be meromorphically 
continued 
to the half plane $\Re(s) >0$.
If $h(X)$ is a product of cyclotomic polynomials, i.e.\ if $|\alpha_j|=1 $ for
every $j$, then and only then can $Z(h;s)$
be continued to the whole complex plane. In all other cases the imaginary axis
is the natural boundary.
\end {Theo}
 The strategy of his proof was to show that every point on the line
$\Re\,s=0$ is an accumulation point of poles or zeros of $\Z$. Estermann's
method was subsequently generalised by many authors.

Dahlquist \cite{Dahl}, for example, extended the above case to $h$ being any
analytic function with isolated 
singularities
within the unit circle. He used the concept of vertex numbers and showed that
except for the case where $h(p^{-s})$ has a finite number of factors of the form
$(1-p^{-\nu s})^{-\beta_{\nu}}$, there is a natural boundary of the zeta
function
at $\Re\,s=0$. 

Later, Kurokawa \cite{Ku} continued on the idea of Estermann to cases where $h$ depends
on the traces of representations of a topological group and solved Linnik's problem for the analytic
continuation of scalar products 
of the
Hecke-$L$ series 
$L(s;\chi_i)$
 where $\chi_i$ are  Gr\"osssencharakters (not necessarily of finite order) 
of finite extensions of an algebraic number field.
His result can be stated more precisely as~:

Let $F/{\Q}$ be a finite extension and $K_i/F$ be $r$ finite extensions of
degree $n_i$ each. The scalar product $L(s;\chi_1,\ldots,\chi_r)$  has the imaginary axis as the
natural boundary except when
\begin{eqnarray*}
(n_1,\ldots,n_r)
&=& (1,\ldots\ldots\ ,1, \star)\qquad \text{or}
 \\
&=&  (1,\ldots .\ ,1,2,2), 
\end{eqnarray*}
in which case  $L(s;\chi_1,\ldots,\chi_r)$ can be continued to the whole
of $\C$ (ibid. Part II, Theorem 4). 

There is of course, no reason to believe that the natural boundary would
always be a line. 
In an example involving the
Euler-phi function  \cite{these} 
$$
Z(s)=\sum_{n=1}^{\infty}\frac{1}{\phi
  (n)^s}=\prod_p\big(1+(p-1)^{-s}(1-p^{-s})^{-1}\big),$$
the boundary of continuation is an open, simply connected, dense set
of the half-plane $\Re\,s>-1$.

The question of analytic continuation of Euler products in several variables 
occur naturally in very many contexts. To cite just one  example, in the study of strings over $p-$adic fields \cite{BFOW}, products of 5-point  amplitudes 
for the  open strings are considered, where the amplitudes are defined as $p-$adic integrals
$$
A_5^p(k_i)=\int_{\Q_p^2}\mid x\mid^{k_1k_2}\mid y\mid^{k_1k_3}\mid 1-x\mid^{k_2k_4}
\mid 1-y\mid^{k_3k_4}\mid x-y\mid^{k_2k_3} dxdy.  $$ 
The product $\prod_p A_5^p$ can be analytically continued to the whole of $\C$, which gives 
interesting relations of such amplitudes with real ones.

We would thus like a multivariable Estermann type of theorem. For this 
we need some notation \cite {BEL}. Let us consider $n$-variable integer polynomials $h_k$
and let 
$$
h(X_1,\dots, X_n, X_{n+1})=1+\sum_{k=0}^d h_k(X_1,\dots, X_n)X_{n+1}^k.
$$
The exponents of the monomial occurring in this expression determine a
polyhedron in $\R^n$ and enable us to give a  description of the domain of
convergence
for the Euler product in $n$ complex variables $Z(h;s_1,\dots,s_n)=
\prod_{p \ \text{prime}} h(p^{-s_1},\dots,p^{-s_n})$. Thus we define, for $\delta \in \R$,  
$$
V(h;\delta):=\bigcap_{k=0}^d\{
    s \in \C^n ~~|~~ \Re(\langle \alpha, s \rangle) >k+\delta  \quad \forall\ 
\alpha \in Ext(h_k)\}
$$ 
where $ Ext(h_k)$ is the set of those points which do not belong to the
interior of any closed segment of the Newton polyhedron
of $h_k$. 
We show that the geometry of the natural boundary is that of a tube over a
convex set with piecewise linear boundary and give a criterion for its
existence which is analogous to Theorem~\ref {est}.

A polynomial $h$ in several variables is called  cyclotomic, if there exists
a finite set  of non-negative integers 
$m_{i,j}$ and a finite set of integers
integers $(\gamma_j)_{j=1,\dots,q}$ such that: 
$$h(X)=\prod_{j=1}^q  
(1-X_1^{m_{1,j}}\dots X_n^{m_{n,j}})^{\gamma_j}.$$
In \cite{BEL} we prove that either $h$ is cyclotomic, or it determines
a natural
boundary of meromorphy, i.e.
\begin{Theo}\label{ester+}

The Euler product
$$ Z(h;s)=Z(h;s_1,\dots,s_n)=
\prod_{p} h(p^{-s_1},\dots,p^{-s_n})
$$
converges absolutely in the domain $V(h;1)$ and can be
meromorphically continued to the domain $V(h;0)$.\par 
Moreover, $Z(h;s)$ can be continued to the whole complex space $\C^n$
if and only if $h$ is cyclotomic.
In all other cases $V(h;0)$ is a natural boundary.

\end{Theo}
Using Newton polyhedra we can write the above as a product of Riemann zeta
functions and a holomorphic function in $V(h;1/r)$, for every natural number
$r$, i.e.
$$
Z(h;s)=\big(\prod_{1\le |m|\le N_r}\zeta(\langle m,s\rangle )^{\gamma(m)}\big)G_{1/r}(s)$$
where $m$ is a $n$-tuple of positive integers, $\{N_r\}$ an increasing
sequence of positive integers and $G(s)$  an absolutely convergent Euler product.
We then treat separately the cases where the set $\{m:\gamma(m)\ne 0\}$ is
finite or infinite
to show that a meromorphic continuation to  $V(h;\delta)$ is not possible for
any $\delta<0$.

A result similar to the above theorem can also be obtained for Euler products of analytic functions
on the unit poly-disc $P(1)$ in $\C^n$
rather than polynomials (op.cit. Theorem 4).
However Theorem~\ref{ester+} is in general not enough to treat Euler products 
of the form $\prod_ph(p, p^{-s})$ which occur, for example, in zeta functions
of groups
and height zeta functions. In certain cases authors have been able to find
natural boundaries of such Euler products while even for an
apparently simple case 
like $f(s) = \prod_p \Big(1+p^{-s} + p^{1-2s}\Big)$ \cite{duSPre}
we might be unable to provide a complete answer (see the next section). 

In fact it does not suffice to
prove that each point is a limit point of poles or zeros of the single factors,
since poles and zeros could cancel. 
In certain situations it is possible to find conditions which ensure that too much 
cancellation among potential singularities is impossible and thereby 
get information on
series like the one just cited.  
For instance, in \cite {natbd} we obtain :
\begin{Theo}
\label{Thm:SpecialBoundary}
Assume the Riemann $\zeta$-function has infinitely many zeros off the line
$\frac{1}{2}+it$. Suppose that $f$ is a function of the form
$f(s)=\prod_{\nu\geq1}\zeta(\nu(s-\frac{1}{2})+\frac{1}{2})^{n_\nu}$ where the
exponents $n_\nu$ are rational integers and the series
$\sum\frac{n_\nu}{2^{\epsilon\nu}}$ converges absolutely for every
$\epsilon>0$. Then $f$ is holomorphic in the 
half plane $\Re s>1$ and has meromorphic continuation in the half plane $\Re
s>\frac{1}{2}$. Denote by $\mathcal P$ the set of prime numbers $p$, such that
$n_p>0$, and suppose that for all $\epsilon>0$ we have
$\mathcal P((1+\epsilon)x)-\mathcal P(x)\gg x^\frac{\sqrt{5}-1}{2}\log^2 x$. Then the line $\Im
s=\frac{1}{2}$ is the natural boundary of $f$; more precisely, every point of
this line is accumulation point of zeros of $f$.
\end{Theo}
To get this result we need some combinatorial geometry on the lines of
Dahlquist \cite{Dahl}.
The following is a sketch of the argument to get the above natural boundary.
By assumption of the falsity of the Riemann hypothesis,
 for every $\epsilon>0$ and every $t$ there is a zero $\rho=\sigma+iT$ of
$\zeta$, such that $\mathcal{P}(T/t)-\mathcal{P}(T/((1+\epsilon)t))\gg
(T/t)^\theta\log^2 (T/t)$, where $\theta=\frac{\sqrt 5-1}2$.
Instead of showing that this particular $\rho$ cannot be cancelled out by poles or zeros of
other factors,, we show that not all zeros can be cancelled out. If
$\frac{\rho-1/2}{p}+\frac{1}{2}$ is not a zero of $f$ for any 
$p\in\mathcal{P}$ and any  $\frac{T}{p}\in[t,
(1+\epsilon)t]$, using combinatorial arguments we reach the contradictory
conclusion that $\theta<\frac{\sqrt{5}-1}{2}$. So in every
square of the form $\{s:\Re\;s\in[\frac{1}{2}, \frac{1}{2}+\epsilon], 
\Im\;s\in[t, t+\epsilon]\}$, there is a zero of $f$.
 
Concerning general Euler products of polynomials in $p$ and  $p^{-s}$, there
exists a conjecture \cite{book}.
\begin{conj}
Let $W(x, y)=\sum_{n,m} a_{n,m} x^n y^m$ be an integral polynomial
with $W(x, 0) = 1$. Then $D(s)=\prod_p W(p, p^{-s})$ is
meromorphically continuable to the whole complex plane if and if only
if it is a finite product of Riemann $\zeta$-functions. Moreover, in
the latter case if $\beta=\max\{\frac{n}{m}: a_{n,m}\neq 0\}$, then
$\Re\,s=\beta$ is the natural boundary of $D$.
\end{conj}
Though all known examples confirm this it is still far from
being resolved. 
In fact we believe that
any refinement of Estermann's method is not enough
to prove this conjecture \cite{aby2}.

We define an obstructing point $z$ to be a complex number with
$\Re\,z=\beta$, such that there exists a sequence of complex numbers
$z_i$, $\Re\,z_i>\beta$, $z_i\rightarrow z$, such that $D$ has a pole
or a zero in $z_i$ for all $i$. Obviously, each obstructing point is
an essential singularity for $D$, the converse not being true in
general. 

Since $D$ may not be convergent on the half-plane
$\Re\,s>\beta$, to continue it meromorphically  
it is written as a product of Riemann $\zeta$-functions and a function $R(s)$
holomorphic, zero-free, and bounded on every half-plane
$\Re\,s>\beta+\epsilon$. Thus  there exist integers
$c_{n,m}$ such that
\[
D(s) = \prod_{n,m} \zeta(ns+m)^{c_{n,m}}\times R(s).
\]
When approximating $D(s)$ by a product of Riemann $\zeta$-functions,
the main contribution comes from monomials $a_{n,m} x^n y^m$ with
$\frac{n}{m}=\beta$. We collect these monomials together in
$\tilde{W}$ that is, we have
\[
W(x, y)=\tilde{W}(x, y) + \sum_{n,m}\nolimits^* a_{n,m} x^n y^m,
\]
where $\sum^*$ means summation over all pairs $n,m$ with
$\frac{n}{m}<\beta$ (in \cite{duSG} the terminology  `ghost polynomial' is used). 

We can classify such polynomials into exclusive, non-empty cases as follows :

\begin{enumerate}
\item $W=\tilde{W}$ and $W$ is cyclotomic; in this case, $D$ is a
  finite product of Riemann $\zeta$-functions;
\item $\tilde{W}$ is not cyclotomic; in this case, every point of the
  line $\Re\,s =\beta$ is an obstruction point;
\item $W\neq\tilde{W}$, $\tilde{W}$ is cyclotomic, and there are
  infinitely many pairs $n,m$ with $a_{n,m}\neq 0$ and
  $\frac{n}{m}<\beta<\frac{n+1}{m}$; in this case, $\beta$ is an
  obstruction point;
\item $W\neq\tilde{W}$, $\tilde{W}$ is cyclotomic, there are only
  finitely many pairs $n,m$ with $a_{n,m}\neq 0$ and
  $\frac{n}{m}<\beta<\frac{n+1}{m}$, but there are infinitely many
  primes $p$ such that the equation $W(p, p^{-s})=0$ has a solution
  $s_0$ with $\Re\,s_0>\beta$; in this case every point of the line
  $\Re\,s=\beta$ is an obstruction point;
\item None of the above; in this case, no point on the line
  $\Re\,s=\beta$ is an obstruction point.
\end{enumerate}

In the third case we need an understanding of the zeros of the
Riemann-zeta function to have information about the meromorphic continuation
and as we will see in the next section that this may only
give conditional answers.  However in
the last case we might be able to say nothing about the analytic continuation
as we will see in the example of $D(s)=\prod_p (1-p^{2-s}+p^{-s})$. 

We would need some
really new
ideas to understand Euler products of polynomials in $p$ and $p^{-s}$.
\subsection{A random series}
From a probabilistic point of view, it is usual to study random
Dirichlet series 
 and show that  almost surely they
have natural
boundaries. Such generic conditions comfort us in 
the belief that for a Dirichlet series there should be meromorphic continuation up to an expected
domain.
 
Often in the definition of a random series 
the coefficients are random (for example
in Kahane \cite{Ka} or Qu\'effelec \cite {Queff}). In the following
\cite{natbd} we use random variables in the exponent to resemble the
Euler
products $W(p,p^{-s})$ discussed before. 

We call a function
regular in a domain if it is meromorphic up to a
discrete set of branch points in the domain,
 that is, it is holomorphic with
the exception of poles and branch points. We can now state the
following probabilistic result :
\begin{Theo}
\label{thm:RandomCont}
Let $(a_\nu), (b_\nu), (c_\nu)$ be real sequences, such that $a_\nu,
b_\nu\to\infty$, and set $\sigma_h=\limsup\limits_{\nu\to\infty}
-\frac{b_\nu}{a_\nu}$. Let $\epsilon_\nu$ be a sequence of independent
real random 
variables, such that

\[
\liminf_{\nu\to\infty} \max_{x\in\R} P(\epsilon_\nu=x) = 0,
\]
and suppose that for $\sigma>\sigma_h$ the series
\begin{equation}
\nonumber\label{eq:ConvCond}
\sum_{\nu=1}^\infty \frac{|c_\nu+\epsilon_\nu|}{2^{a_\nu\sigma+b_\nu}}
\end{equation}
converges almost surely. Then with probability 1 the function
\[
Z(s) = \prod\limits_{\nu=1}^\infty \zeta(a_\nu s+b_\nu)^{c_\nu+\epsilon_\nu}
\]
is regular
in the half-plane $\Re\;s>\sigma_h$ and has the line
$\Re\;s=\sigma_h$ as its natural boundary.
\end{Theo}
To give an idea of the arguments used in the proof we let $s_0=\sigma_h+it$ be a point on the supposed
boundary with $t\neq 0$ rational, and consider the
square $S$
with side length $\frac{2}{n}$ centred in $s_0$. For $\epsilon>0$  given, we show that with probability
$>1-\epsilon$ the function $Z$ is either not meromorphic on $S$, or has a zero or a
pole in $S$. 
 Then for
a suitably chosen index $\mu$ we consider\[
Z_\mu(s) = \prod\limits_{\nu\neq\mu}^\infty \zeta(a_\nu
s+b_\nu)^{c_\nu+\epsilon_\nu}.
\]
such that if $Z$ is meromorphic on $S$, so is $Z_\mu$. Let $D_1$ be the divisor of the
restriction of $Z_\mu$ to $S$, and let $D_2$ be the divisor of $\zeta(a_\mu
s+b_\mu)$ restricted to $S$. We  show that
$D_1+(c_\mu+\epsilon_\mu)D_2$ is non-trivial with probability
$>1-\epsilon$. 
 The number
of zeros of $\zeta(a_\mu s+b_\mu)$ in $S$ equals $N(T+h)-N(T)$, where $N$
denotes the number of zeros of $\zeta$ with imaginary part $\leq T$, and $T$
and $h$ are certain real numbers satisfying $T\geq 1000$ and $h\geq
6$. Using a classical estimate \cite{Ba}, we can show that $D_2$ is non-trivial.

We note that in the initial statement of the above theorem, the term
 `holomorphic' appeared instead of  `regular' (Theorem 3, ibid.). 
In fact, as pointed out by J-P. Kahane, the finite product of
$\zeta$-functions dominating the behaviour of $Z$ in a half-plane
$\Re\;s>\sigma_h+\epsilon$ can yield branch points at all poles and
zeros of the involved $\zeta$-functions.

\section{Examples}
\subsection{Zeta function of a symplectic group}
The local
zeta function associated to the algebraic group $\G$  is defined as 
$$
Z_p(\G, s)=\int_{\G_p^+} \mid  \det(g)\mid_p^{-s}d\mu
$$
where $\G_p^+=\G(\Q_p)\cap M_n (\Z_p)$ , $\mid.\mid_p$ denotes the 
p-adic valuation
and $\mu$ is the normalised Haar measure on $ \G(\Z_p)$.
In
\cite{duSG}, du Sautoy and Grunewald prove that the natural boundary of the
zeta function $Z(\G,s)$ of the symplectic group $GSp_6$  given by \cite{Ig}
\begin{eqnarray*}
Z(s/3) & = & \zeta(s)\zeta(s-3)\zeta(s-5)\zeta(s-6)\prod_p
\Big(1+p^{1-s}+p^{2-s}+p^{3-s}+p^{4-s}+p^{5-2s}\Big)\\
\end{eqnarray*}
has a natural boundary at $\Re\;s=\frac 4{3}$.
To show that every point on the boundary is an accumulation point of zeros, 
the authors consider the partial derivatives at $(0,-1)$ of the equation
$1+(1+V+V^2+V^3)U+V^3U^2=0$ where $ V=p^{-1},\ U=p^{4-s}$. The Implicit
Function Theorem then guarantees the existence of a solution  for the above
equation in $U=-1+V+\Omega (V)$, where for $p$ large enough $\Omega (V)$ contains terms smaller than
$p^{-2}$. Thus for every integer $n$ there is a solution 
$$s=4-\frac{\log(1-p^{-1}+\Omega(p^{-1}))}{\log p}+\frac{(2n-1)\pi i} {\log p}.$$
Now for a large prime $p$ and a fixed point $A$ with $\Re\;s=4$ on the boundary 
we can find a sequence of integers $n_p$
such that $\frac{(2n_p-1)\pi} {\log p}\rightarrow \Im (A)$. Further the fact
that 
$$-\frac{\log(1-p^{-1}+\Omega(p^{-1}))}{\log p}>0$$ 
for large enough $p$
means that  $Z(s/3)$ cannot be continued beyond its assumed boundary
$\Re\;s= 4$.

Notice that this is an example of the  `lucky' situation we
encountered in the fourth case of the classification of $\prod W(p,p^{-s})$.

\subsection{Height zeta functions}

Several people in the recent past have studied the analytic
properties of 
height zeta functions associated 
to counting rational points on algebraic varieties. Of particular interest is
the case of  a variety with ample anti-canonical bundle (called a Fano
variety) $V$
over a number field $k$ whose $k$-rational points are Zariski dense in $V$,
for a  height function $H$ defined naturally over the anti-canonical sheaf. 
Here  an important motivation is Manin's conjecture
that, for $U$ a suitably defined open subset of $V$, 
$$ |\{x\in U(k) :H(x)\le t\}|\sim Ct(\log t)^{r-1} $$
 as  $t\rightarrow\infty$.
In the above, $C$ is a non zero constant and $r$ the rank of the Picard group of $V$.
There is a further conjecture, due to Peyre \cite{Peyre}, on the constant $C$ relating it to the Tamagawa measure.

We will concentrate on  
 $\P^n$, the projective $n$-space over the field $\Q$ with the
classical normalised height function $H_n : \P^{n-1}\rightarrow \R_{>0}$ defined by 
$H(x )= \max_i\{|x_i|\}$, where $\gcd (x_1,\ldots ,x_{n+1})=1$. (Other definitions of the height exist but
we shall not treat them here. The interested reader could see, for
example, \cite{Ess}). 

We now give details of analytic continuation and boundaries of
a few zeta functions in the above context which have Euler products in 
several variables.

\subsubsection {A cubic surface}  For studying the above case, it is possible to first choose the anti-canonical
line bundle and assume that it be ample. This then determines a projective
embedding of the desingularised model of the variety using a fan
decomposition into finitely many simplical integral cones (for details see, for example,
\cite{bret1} or \cite{Salb}).
The zeta function is then defined, for $\Re\,s$ large enough, as 
$$Z_{U}(s)=\sum_{x\in U}\frac{1}{H(x)^{s}}.$$
De la Breteche and
Swinnerton-Dyer \cite{2nobles} proved that the zeta function associated to singular cubic
surfaces has a natural boundary at $\Re\,s=3/4$.  We follow the treatment of
\cite{bret2} to give a summary of their original proof where
they study the multi-variable
function
$$Z(s_1,s_2,s_3)=\underset{gcd(x_1,x_2,x_3)=1}
{\sum_{x_1x_2x_3=x_4^3}}
   \frac
{1}{x_1^{s_1}x_2^{s_2}x_3^{s_3}}$$
outside the union of three lines in the hyper-surface $x_4=0$.
The  Euler product of this function is given by
$$\prod_p\frac{1+\sum(1-p^{3s_i})p^{-(2s_j+s_k)}-p^{-(3s_1+3s_2+3s_3)}}
{(1-p^{3s_1})(1-p^{3s_2})(1-p^{3s_3})}$$
where in the sum each of $i,j,k$ take the values $1,2,3$.
The above is then written with the help of functions  `convenient' for 
$\Re\, s>\frac{3}{4}$ 
and $F(s)$ which involves the Euler product of a polynomial in two variables
$$W(x,y)=1+(1-x^3y)(x^6y^{-2}+x^5y^{-1}+x^4+x^2y^{2}+xy^{3}+y^{4})-x^9y^3$$
with  $x=p^{-1/4},\  y=p^{3/4-s}$.
Here again the authors succeed in establishing that every point on the assumed boundary
is the limit point of  a subset of zeros of the
function $F(s)=\prod_pW(p^{-1/4}, p^{3/4-s})\prod_{j=-2}^4(1-p^{-(1+j(s-1)})$.
For a fixed prime $p$, the  number of zeros with $\Re\,s>\frac{3}{4}$
of $W$, i.e. 
$$ \frac{3}{4}+ \frac{1}{4\sqrt 2p^{1/4}\log p}+\mathcal O\big(
\frac{1}{p^{3/4}\log p}\big)$$
is large.
Now for $\Re s>\frac{3}{4}+\frac{1}{N}$, 
there exist suitably chosen finite number of integers $b(k,k')$ such that
$$
F(s)=\underset{ k-k'/4+k'/N>1}{\prod_{k,k'}}\zeta(k+k'(s-1))^{b(k,k')}\prod_pW_N(p^{-1/4}, p^{3/4-s})$$
where 
$$
W_N(p^{-1/4}, p^{3/4-s})=W(p^{-1/4}, p^{3/4-s})\underset{ k-k'/4+k'/N>1}{\prod_{k,k'}}
 (1-p^{(k+k'(s-1))})^{b(k,k')}.$$
The zeros of $W_N(p^{-1/4}, p^{3/4-s})$ and $W(p^{-1/4}, p^{3/4-s})$ are the
same. Further for every real $\tau$ one can construct a sub-sequence of its
zeros which converge to $\frac{3}{4}+i\tau$ and which are not poles of 
$$ \underset{k-k'/4+k'/N>1}{\prod_{k,k'}}
 \zeta(k+k'(s-1))^{b(k,k')}.$$
These zeros are again the zeros of $F(s)$ and therefore no continuation is
possible beyond the assumed boundary.

For what concerns the asymptotics, it is known that
$$|\big\{x\in U: H(x)\le t\big\}|=tQ(\log
t)+\mathcal{O}(t^{7/8+\epsilon})$$
where the degree of $Q$ is 6 and the leading coefficient is
$\frac{1}{6}\prod_p\{(1-1/p)^7(1+7/p+1/p²)\}$.
\subsubsection {An n-fold product}

In \cite{BEL}, we consider instead an implicit projective embedding determined by  a
finite set of equations and do not need a fan decomposition.

Let $X$ be a toric variety and $A_{d,n}=A$ a $d\times n$ integer matrix all of
whose row sums are zero.
 The rational points of the toric variety are defined by 

$$X(A):=\{(x_1,\cdots, x_n)\in \P^{n-1}(\Q) : \prod_{i:a_{j,i}\ge
    0}x_i^{a_{j,i}}=\prod_{i:a_{j,i}<0}x_i^{-a_{j,i}}\ \forall j\}$$
and the maximal torus $U(A)$ comprises of those elements of $X(A)$ the product of whose
coordinates is non-zero.
Each point in the maximal torus corresponds to a unique $n$-tuple of co-prime positive
integers which we denote by $(m_1,\cdots,m_n)$.
 
We define a multivariable zeta function, for $\Re\,s_i>1 $, comparable to the one used 
for toric varieties in \cite{bret1} as
$$ Z_A(s)=\sum_{m_i\in \N} \frac {F_A(m_1,\cdots,m_n)}{m_1^{s_1}\dots
  m_n^{s_n}},$$
where 
\begin{eqnarray*}
F_A(m_1,\cdots,m_n)
&=& 1\qquad\text{if}\ \gcd\ (m_1,\cdots,m_n)=1 ,\ \prod_i
m_i^{a_{j,i}}=1\ \forall j,
\\
&=&  0\qquad\text{otherwise.}
\end{eqnarray*}

The defining equations are multiplicative and we thus get an Euler product
expansion of an analytic function in $n$ complex variables
$$ Z_A(s)=\prod_p h_A(p^{s_1},\dots,p^{s_n}).$$
Now the function $h_A(X)$ is expressed as a rational function
$$\prod_{\nu \in K}(1-X^{\nu})^{-c(\nu)}W(X)$$
for positive integers $c(\nu)$, a finite index set $K$ and an integer $n$-variable polynomial $W$.
We can prove, using Theorem~\ref{ester+}, that $Z_A(s)$ has a natural
boundary. 
In fact, it is possible to explicitly describe the whole boundary of analytic
continuation
(see \cite{BEL},  Theorem 6).
The description of the analytic continuation of this zeta function can now be used to deduce the
asymptotic properties of the height density function on $U(A)$ because of the equation
$$ |\big\{x\in U(A): H(x)=\max_i|m_i(x)|\le t\big\}|=C(A)\sum_{m_i\le t}F_A(m_1,\cdots,m_n)$$
where $C(A)$ is a computable constant.

As a special case, we get asymptotic results for the number
of $n$-fold products of relatively prime positive integers that equal the 
$n$th power of an integer. Batyrev and Tschinkel \cite{BT} showed that this
problem is equivalent to the 
asymptotic description of the height
density function on the maximal torus of the hyper-surface $ x_1\cdots x_n=x_{n+1}^n$.
Now that there is only one equation involved, i.e. $d=1$, we use the
matrix $A_n=(1,\cdots,1,-n)$ and the rational points are
$$U(A_n):=\{(x_1,\cdots, x_{n+1})\in \P^{n}(\Q) : x_1\cdots x_n=x_{n+1}^n,,\
x_1\cdots x_n\ne 0\}.$$
To express $h_{A_n}(X)$ precisely as a rational function on the unit poly-disc
$P(1)$  we notice
that if $$h_{A_n}(X)=\sum_{\alpha}X^{\alpha} $$ 
for $\alpha\in \N_0^{n+1}$, to satisfy the
definition of $F_{A_n}$ we require that $A_n(\alpha)=0$. Thus  
for a $n$-tuple $r$, we use the notation $|r|$
for its weight, i.e.\ the sum of its
$n$ coordinates and for all $r$ such that $|r|/n$ is a non-negative
integer, we let  
$l(r)=(r_1,\ldots,r_n,|r|/n)$. Further we ensure that 
the condition of coprimality of the $p^{\alpha_i}$ is met and obtain
$$h_{A_n}(X)=(\prod_{i=1}^n(1-X_i^nX_{n+1})^{-1})\sum_{|r|/n\in \N_0}
  X_1^{r_1}\ldots X_n^{r_n}X_{n+1}^{|r|/n}.$$
The sum in the expression above is not cyclotomic and this gives the
natural boundary of $$Z_{A_n}(s)=R\times\prod_p(\sum_{r \in D_n}\frac{1}{p^{\langle
    l(r), s \rangle}})$$
for $R$ a finite product of Riemann zeta functions, to be
$$ V(0)=\{
    s \in \C^{n+1} ~~|~~ \Re(\langle l(r), s \rangle) >0  \quad \forall\ 
r \in D_n\}.
$$ 
Using the above analytic properties and a multivariable
 Tauberian theorem \cite{bret1} we prove
that
\begin{Theo} 
There exists $\theta>0$ such that
$$|\big\{x\in U(A_n): H(x)\le t\big\}|=tQ_n(\log t)+\mathcal{O}(t^{1-\theta})$$
where $Q_n(\log t)$ is a non-vanishing polynomial of degree $d_n=\binom{2n-1}{n}-n-1$.
\end{Theo}
Actually we can describe the last polynomial rather precisely
for all $n\ge 3$ (\cite{BEL}, Theorem 7).
\subsection{Unlucky cases}

In the last two subsections we could give satisfactory descriptions of
the analytic behaviour of the Euler products. This need not always be
possible.
In the following we can only show the existence of a conditional natural
boundary \cite{natbd}.

\begin{Prop}\label{prop:innocent1}
Suppose that there are infinitely many zeros of $\zeta$ off the line
$\frac{1}{2}+it$. Then the function 
\[
f(s) = \prod_p \Big(1+p^{-s} + p^{1-2s}\Big)
\]
has meromorphic continuation to the half plane $\Re s>\frac{1}{2}$, and the
line $\Re s=\frac{1}{2}$ is the natural boundary of $f$. 
\end{Prop}
This is an example of case (3) of our classification of the previous
section. We notice that the real parts of the zeros of 
$f(s)$
are exactly $\frac{1}{2}$
and thus we can not construct a sub-sequence of zeros or poles  which
would converge on each point of the presumed natural boundary $\Re
s=\frac{1}{2}$. 
The conditional result above is attained by expressing $f(s)$ as a
product of 
functions  `convenient' for 
$\Re s>\frac{1}{2}$ 
and $$\prod_{m\ge 1} \frac{\zeta((4m+1)s-2m)}{\zeta((4m+3)s-2m-1)},$$ which
is of the type considered in Theorem~\ref{Thm:SpecialBoundary}.

We next
consider the 
Euler product  
$D(s)=\prod_p (1-p^{2-s}+p^{-s})$ 
which can be written as
\[
D(s) = \prod_p (1-p^{2-s}) \prod_p (1+\frac{p^{-s}}{1-p^{2-s}}) 
= \zeta(s-2)D^*(s).
\]
We expect a natural boundary at $\Re\;s=2$, or at least an
essential singularity at $s=2$, and our only method to prove this is
to approach this point from the right. But for $\Re\;s=\sigma>2$ we estimate
the second  product as 
\[
\sum_p \sum_{n\geq 1} |p^{2n-(2n+1)s}| \leq \sum_p p^{-2} \sum_{n\geq
  1} p^{-n(\sigma-2)} \leq \sum_p p^{-2} \frac{1}{1-2^{\sigma-2}}.
\]
So the product for $D^*$ converges absolutely in the half-plane
$\Re\;s>2$, in particular, $D^*$ does not have any zeros or poles in
this half-plane.
This example falls under case (5) of the classification mentioned.
It is worse than  Proposition~\ref{prop:innocent1} where we
could not unconditionally prove  the existence of zeros or poles clustering on
the assumed boundary  whereas here such zeros or poles
do not even exist.

\section{No Euler products}
There are numerous contexts in which we come across zeta functions that do not
have an Euler product.
We cite just two examples. The first, mentioned because it comes from a context
quite different from the other examples we treated, is that of Dirichlet series generated by
finite automata.
  
Roughly speaking, a sequence $(u_n)$ with values in a finite set is
$d$-automatic 
if we can compute the $n$-th
term of the sequence by feeding
the base $d$ representation of $n$ to a finite state machine. One of the best
known  among $2$-automatic cases is 
the Thue-Morse sequence, 
$$ 01\ 10\ 1001\ 10010110\ \cdots $$
generated by the substitution maps $0\rightarrow 01,\
1\rightarrow 10$. The Dirichlet series
$\sum_{n=0}^\infty \frac {u_n}{n^s}$ corresponding to a $d$-automatic
sequence can be meromorphically continued to the 
whole complex plane. Among consequences it is proved  \cite{All} that automatic sequences have
logarithmic densities. It would be interesting to know how Dirichlet
series associated to non automatic sequences (like the infinite Fibonacci word
generated by the substitutions $0\rightarrow 01,\
1\rightarrow 0$) behave.

The second example is in several variables. 
The Euler-Zagier sum defined as 
$$ \zeta_r(s_1,\cdots,s_r)=\sum_{m_1=1}^{\infty}\cdots\sum_{m_r=1}^{\infty}
m_1^{-s_1}(m_1+m_2)^{-s_2}\cdots(m_1+\cdots+m_r)^{-s_r}$$
has been studied with much enthusiasm. This function can be analytically
continued to the whole $\C^r$ space.
Matsumoto  introduced the generalised multiple zeta function 
$$ \zeta_r((s_1,\cdots,s_r);(\alpha_1,\cdots,\alpha_r),(w_1,\cdots,w_r))=
\sum_{j=1}^r\sum_{m_j=0}^{\infty}\prod_{i=1}^r
(\alpha_i+m_1w_1\cdots+m_iw_i)^{-s_i}
$$
where $w_i, m_i$ are complex parameters with branches of logarithms suitably defined.
This too can be continued as a meromorphic function to the whole $\C^r$ plane.
We do not wish to elaborate on this subject but the interested reader can find
details  elsewhere (see, for example, the 
expository paper \cite{mat} for references).

\subsection {Goldbach zeta function}

Here we consider the number $G_r(n),\ r\ge 2 $, of
representations of $n$ as the sum of $r$ primes. 
Egami and
Matsumoto\cite{EgMat} introduced the generating function
\[
\Phi_r(s) = \sum_{k_1=1}^\infty\dots\sum_{k_r=1}^ \infty 
\frac{\Lambda(k_1)\dots\Lambda(k_r)}{(k_1+k_2+\dots+k_r)^s} 
 = \sum_{n=1}^\infty\frac{G_r(n)}{n^s}
\] 
using the von
Mangoldt function $\Lambda$.
This series is absolutely convergent for $\Re\;s>r$, and has a simple pole
at $s=r$. It is clear that to study the analytic properties in this 
context it is necessary to have information on the zero-free region of
$\zeta$, the
Riemann zeta function, and the presence of even one zero of $\zeta$ may
prevent us from having useful information. All results that we will talk
about will therefore be under the assumption of the Riemann Hypothesis (RH).

We can show that from the analytic point of view, under RH, $\Phi_r$  is
determined by the case $r=2$
\cite{gold}.
\begin{Theo}
\label{thm:Rep}
Suppose that the Riemann Hypothesis is true. Then for any $r\geq 3$ 
there exist polynomials $f_r(s), g_r(s),
h_r(s)$, such that
\[
\Phi_r(s) = f_r(s)\zeta(s-r+1) + g_r(s)\frac{\zeta'}{\zeta}(s-r+1) +
h_r(s)\Phi_2(s-r+2) + R(s),
\]
where $R(s)$ is holomorphic in the half-plane $\Re s>r-1$ and
uniformly bounded in each half-strip of the form $\Re s>r+1$, $T<\Im
s<T+1$, with $T>0$.
\end{Theo}
This is done by computing the function using the circle method which give the
three
main terms.
A bound (under RH) for
 $$\sum_{n\leq x} \Lambda(n)e^{2\pi i\alpha n}-
\sum_{n\leq x} e^{2\pi i\alpha n}$$
gives an error term of order $\mathcal O(x^{r-1-\delta})$ for some $\delta$
positive for all but the above three terms. 

It is thus important to consider the situation of $r=2$. We recall that this
case occurs 
in the consideration of the Goldbach conjecture
that every  even integer larger than $2$ is the sum of two primes. To
study this problem often
it is natural to consider the corresponding problem for $\Lambda$ and
try to
show that $G_2(n)>C\sqrt n$.

Now,
assuming the RH, the authors in
\cite{EgMat} described the analytic continuation of $\Phi_2$
and  
for obtaining a natural boundary they used unproved assumptions on the
distribution of the imaginary parts of zeros of $\zeta$. 
In this context we denote  the set of imaginary parts of non-trivial
zeros of $\zeta$ by $\Gamma$. The belief that the positive elements in
$\Gamma$ are rationally independent is folkloric and
Fujii\cite{Fu1} used the following special case :
\begin{conj}
\label{Con:Indep}
Suppose that $\gamma_1+\gamma_2=\gamma_3+\gamma_4\neq 0$ with
$\gamma_i\in\Gamma$. Then $\{\gamma_1, \gamma_2\}=\{\gamma_3,
\gamma_4\}$.
\end{conj}
In \cite{EgMat} an effective version of the above conjecture is
formulated, i.e. 
 \begin{conj} 
\label{Con:IndepEff}
There is some $\alpha<\frac{\pi}{2}$, such that for $\gamma_1, \ldots,
\gamma_4\in\Gamma$ we have either $\{\gamma_1,
\gamma_2\}=\{\gamma_3,\gamma_4\}$, or 
\[
|(\gamma_1+\gamma_2)-(\gamma_3+\gamma_4)|\geq\exp\big(-\alpha
(|\gamma_1| + |\gamma_2| + |\gamma_3| + |\gamma_4|)\big),
\]
\end{conj}
and  it is proven that:
\begin{Theo}
\label{contgold}
Suppose the Riemann hypothesis holds true. Then $\Phi_2(s)$ can be
meromorphically continued into the half-plane $\Re s>1$ with an
infinitude of poles on the line $\frac{3}{2}+it$. If in
addition Conjecture~\ref{Con:IndepEff} holds true, then the line
$\Re\;s=1$ is the natural boundary of $\Phi_2$. More precisely, the
set of points $1+i\kappa$ with $\lim_{\sigma\searrow 1}
|\Phi(\sigma+\kappa)| =\infty$ is dense on $\R$.
\end{Theo}

In [ibid.] the authors conjectured that under the same
assumptions the domain of meromorphic continuation of $\Phi_r$ should
be the half-plane $\Re s>r-1$. Notice that a direct consequence of 
Theorem~\ref{thm:Rep} confirms the following :
\begin{Theo} 
If the RH holds true, then $\Phi_r(s)$ has a natural boundary at
  $\Re s=r-1$ for all $r\geq 2$ if and only if $\Phi_2(s)$ has a
  natural boundary at $\Re s=1$.
\end{Theo}
In \cite{gold} it is also shown that if the RH and
 Conjecture~\ref{Con:Indep} 
hold true, then
  $\Phi_2(s)$ does have a natural boundary at $\Re s=1$ and 
a singularity can be described precisely as

\begin{Theo}
\label{thm:Boundary}
If the RH holds true, then $\Phi_2$ has a singularity at $2\rho_1$, where
$\rho_1=\frac{1}{2}+14.1347\ldots i$ is the first root of $\zeta$. Moreover,
\begin{equation}\nonumber
\lim_{\sigma_\searrow 0} (\sigma-1)|\phi_2(2\rho_1+\sigma)|>0.
\end{equation}
\end{Theo}
 This last result helps us obtain an $\Omega$- result for $G_r(n)$. We
 consider the oscillating term
$$H_r(x)=-r\sum_{\rho}\frac{x^{r-1+\rho}}{\rho(1+\rho)\dots (r-1+\rho)}$$
 where the summation  
runs over all non-trivial zeros of $\zeta$. The generating Dirichlet
series for $$\sum_{n\leq x} G_r(n) -\frac{1}{r!}x^r -H_r(x)$$ has a
singularity at $2\rho_1 +r-1$, which gives the following :
\begin{Cor}
Suppose that RH holds true. Then we have
\[
\sum_{n\leq x} G_r(n) = \frac{1}{r!}x^r + H_r(x) +
\Omega(x^{r-1}).
\]
\end{Cor}

In fact 
the quality of the error term does not improve with  increasing $r$.
We mention a few historical facts about the error term.
Fujii\cite{Fu1} obtained under the RH
$$\sum_{n\le x}G_2(n)=x^2/2+\mathcal{O}(x^{3/2})$$
which he later improved \cite{Fu2}, by explicitly writing the
oscillating term, to 
$$
\sum_{n\le x}G_2(n)=x^2/2+H_2(x)+\mathcal{O}((x\log x)^{4/3}).
$$

Further, in \cite{Nagoya} we used the
distribution of primes in short intervals to estimate exponential sums
close to the point 0 and  proved that
\begin{Theo}
\label{thm:ErrorG2}
Suppose that the RH is true. Then
we have 
\[
\sum_{n\leq x} G_2(n) = \frac{1}{2}x^{2} + H_2(x) +
\mathcal{O}(x\log^5 x),
\]
and
\[
\sum_{n\leq x} G_2(n) = \frac{1}{2}x^{2} + H_2(x) + \Omega_+ (x\log\log x).
\]
\end{Theo}
which confirms the conjectural value of the error term
\cite[Conj.~2.2]{EgMat}.
 
Recently Granville\cite{Gran} used
the error term $\mathcal{O}((x\log x)^{4/3})$ to obtain a new
characterisation of the RH, i.e.\
for the  `twin prime constant' $C_2$, the Riemann Hypothesis is equivalent to the estimate
$$\sum_{n\le
  x}G_2(n)-nC_2\prod_{p|n}\frac{p-1}{p-2}\ll
x^{3/2 +o(1)}.$$ 
 
Using the GRH
one could similarly find bounds for the exponential sums in
question which would give
\begin{Theo}
The Generalised Riemann Hypothesis for Dirichlet $L$-functions $L(s,\chi)$, 
$\chi\bmod{q}$, is 
equivalent to the estimate 
$$\sum_{n\le
  x,\ q|n}G_2(n)=\frac{1}{\phi (q)}\sum_{n\le
  x}G_2(n)+\mathcal O(x^{1+o(1)}),$$ 
\end{Theo} 

as announced in [ibid.Theorem~1C].

\section {Consequences}
One of the amusing consequences of the existence of a natural boundary 
is to suggest that there is a  `natural'
limit to what can be achieved for asymptotic results associated to Dirichlet
series by using complex analysis.
A natural boundary could show 
the non-existence of certain asymptotic results involving error terms
and thus imply the existence of an inverse result, i.e.\ an $\Omega$-term.
Usually when proving an $\Omega$-result we first derive an explicit
formula with oscillating terms and then show that these terms cannot
cancel each other out for all choices of the parameters. In
\cite{natbd} 
we show
that even if we allow for infinite oscillatory sums to be part of the
main terms, we still get lower bounds for the error terms. Thus a
natural boundary at $\Re s=\sigma$ precludes the existence of an
explicit formula with main terms over
the
zeros of the Riemann zeta
function and  an error term $\mathcal{O}(x^\sigma)$. We state this
more precisely as :

\begin{Prop}
\label{NoExplicite}
Let $a_n$ be a sequence of complex numbers such that the generating
Dirichlet-series has a natural boundary at $\Re s=\sigma_h$. Then there does
not exist an explicit formula of the form 
$$A(x) := \sum_{n\leq x} a_n = \sum_{\rho\in\mathcal{R}} c_\rho x^\rho
+\mathcal{O}(x^\theta)$$ 
for any sequence $c$ with $|c_\rho|\ll(1+|\rho|)^c$ and
$|\mathcal{R}\cap\{s:\Re s>\theta,|\Im s|<T\}|\ll T^c$ and for any 
$\theta<\sigma_h$. In
particular, for any sequence $\alpha_i, \beta_i$, $1\leq i\leq k$  and any
$\epsilon >0$ we have
\[
A(x) = \sum \alpha_i x^{\beta_i} + \Omega(x^{\sigma_h -\epsilon}).
\]
\end{Prop}
In practice, the integral taken over the shifted path
 need not always converge and we may not be able
to obtain an explicit formula. This can happen even when the series is meromorphic in the
entire plane, for example, the age-old divisor
problem where we have an $\Omega$-estimate of size $x^{1/4}$ though the
corresponding Dirichlet-series $\zeta^2(s)$ is meromorphic on $\C$. 
 
However, in certain cases
we can actually obtain explicit formulae if
we find bounds on the growth of the Dirichlet-series. We consider a
case
of an Euler product $\prod_pW(p, p^{-s})$ which we have already
encountered as
the $p$-adic zeta function of $GSp_6$. This can be interpreted as a counting
function
\cite{BG} by establishing a bijection between right cosets 
of $2t\times 2t$ symplectic matrices
and sub modules of finite index of $\Z^{2t}$ which are equal to their duals
and called polarised. 

\begin{Theo}

\label{thm:Axexplicite}
Denote by $a_n$ the number of polarised sub modules of
$\Z^6$ of order $n$. Then we have for every $\epsilon>0$ 
\begin{equation}
\label{eq:Axexplicite}
A(x) := \sum_{n\geq 1} a_n e^{-n/x} = c_1 x^{7/3} + c_2 x^2 + c_3 x^{5/3} +
\sum_\rho \alpha_\rho x^{\frac{\rho+8}{6}} + \mathcal{O}(x^{4/3+\epsilon}), 
\end{equation}
where $\rho$ runs over all zeros of $\zeta$, and the coefficients
$c_1$, $c_2$, $c_3$, and $\alpha_\rho$ are numerically computable
constants. Moreover, the error term cannot be improved to
$\mathcal{O}(x^{4/3-\epsilon})$ for any fixed
$\epsilon>0$.
\end{Theo}
The interpretation above allows us to use the zeta function $Z(GSp_6,s)$ 
as the generating function for $a_n$.   
Applying the
Mellin transform
we obtain
\[
A(x) = \frac{1}{2\pi i}\int\limits_{3-i\infty}^{3+i\infty} Z(s)\Gamma(s)x^s\;ds.
\]
For $\sigma$ and $\epsilon>0$ fixed, we have $\Gamma(\sigma+it)\ll
e^{-(\frac{\pi}{2}-\epsilon) t}$. We now choose a path (following
Tur\'an 
\cite[Appendix G]{TuranAna}) to shift the integration. The integral on this new path is bounded above by
$x^{4/3+\epsilon}$. Hence, we obtain the formula
\[
A(x) = \sum_{\Re\rho>4/3+\epsilon} \Gamma(\rho)
x^\rho\res_{s=\rho} Z(s) + \mathcal{O}(x^{4/3+\epsilon}),
\]
where $\rho$ runs over the poles of $Z(s)$, and all complex numbers
$4/3+\rho/6$.
We already saw that $\Re\;s=\frac{4}{3}$ is the natural boundary for
$Z(s)$ and as in Proposition~\ref{NoExplicite}, we now get an $\Omega$-result.

The moral of the story is not necessarily to get the best possible $\Omega$-result
but to show that a non-trivial result is obtainable by this method.

\bigskip
\begin{tabular}{ll}

Universit\'e de Lille 1, \\
Laboratoire Paul Painlev\'e,\\
U.M.R. CNRS 8524, \\
  59655 Villeneuve d'Ascq Cedex, \\
  France.\\
bhowmik@math.univ-lille1.fr
\end{tabular}
\end{document}